\documentclass[12pt]{article}
\usepackage[cp866]{inputenc}
\usepackage[russian]{babel}
\usepackage{amsfonts}
\usepackage{amssymb,amsmath}
\newcommand{\codim}{\operatorname{codim}}
\newcommand{\Cl}{\operatorname{Cl}}
\newcommand{\Si}{\operatorname{Si}}
\newcommand{\dom}{\operatorname{dom}}
\usepackage{amsfonts}\topmargin=-1.5cm\textwidth=17cm \oddsidemargin=-0mm \textheight=24.0cm
\begin{document}

    {\large{\centerline{\bf Asymptotical behavior of
    subspaces under action of}}}
 {\large{\centerline{\bf ~asymptotically finite-dimensional semigroups of operators}}}\vskip3mm

\centerline{ Konstantin Storozhuk}\vskip3mm


Abstract:   {\it We study a~semigroup~$\varphi$
of linear operators acting on a~Banach space~$X$
which satisfies the condition
$\codim X_0<\infty$, where
$X_0=\{x\in X \mid \varphi_t(x)\underset{t\to\infty}\longrightarrow 0\}.$
We show that $X_0$ is closed under these conditions.
We establish some properties concerning the asymptotic behavior
of subspaces which complement~$X_0$ in~$X$.}\vskip3mm

\centerline{\bf 0. Preliminary Definitions and Statements of the Results}\vskip2mm

Let $X$ be a~Banach space and let
$\{\varphi_t:X\to X\mid t\geq 0\}$
be a~semigroup of linear operators; i.e.,
$\varphi_t\circ\varphi_q=\varphi_{t+q}$.
Throughout the article we suppose that the semigroup acts
continuously for
$0<t<\infty$; i.e.,  for each vector $v\in X$ the function
$t\mapsto \varphi_t(v)$ is continuous for $t>0$.
A~semigroup is {\it bounded\/} if all operators $\varphi_t$ are bounded
in the norm by some constant $C<\infty$.

For every vector $v\in X$ we write
$v_t=\varphi_t(v)$
and use the same shorthand notation for arbitrary subsets in~$X$.

Put $X_0=\{x\in X \mid x_t\underset{t\to\infty}\longrightarrow 0\}$.
The space $X_0$ is $\varphi_t$-invariant; i.e.,
$ \varphi_t(X_0)\subset X_0$.

We say that a~semigroup
is {\it asymptotically finite-dimensional\/} if $\codim X_0=n<\infty$.

Although the space $X_0$ consists of vectors tending to zero in the limit,
the restriction of the semigroup to~$X_0$
may fail to be a~bounded semigroup
(and the space $X_0$ may fail to be a~closed subspace in~$X$), see {\bf Example~1}.
However, for asymptotically countably-dimensional
semigroup, $X_0$ is closed and the restriction~
$\varphi_t|_{X_0}$
is a~bounded semigroup even if $\varphi_t$ is not bounded on the whole~$X$
({\bf Theorem~1}). For $X_0$ to be closed, it is also sufficient
that $X$ has a~{\it closed\/} subspace~$Y$ which complements~$X_0$.

Suppose that $Y$ is such a~subspace.
It follows from  the Banach principle that the norm in $X=X_0\oplus Y$
is equivalent to the norm given by the formula $\|(x_0+y)\|:=|x_0|+|y|$.

Since $\varphi_t(X_0)\subset X_0$,
the decomposition of the operators $\varphi_t:X\to X$
has the form
\begin{equation}
\varphi_t=\left(\begin{array}{cc}
\alpha_t & b_t \\
0 & Q_t
\end{array}\right)
:X_0\times Y\to X_0\times Y.
\end{equation}
Existence of a~$\varphi_t$-invariant subspace~$Y$
is obviously equivalent to existence of a~diagonal representation~(1).

The {\it angle} ({\it span}) between two subspaces~$A$
and~$B \subset X$ is the value
$$
\angle(A,B)=\min\{\sup\limits_{a\in A,|a|=1}\{\rho(a,B)\},
\sup\limits_{b\in B,|b|=1}\{\rho(b,A)\}\}.
$$

The angle plays the role of a~metric on
the set of $n$-dimensional subspaces of the space~$X$.

In {\bf Theorem~2} we prove that if $\varphi_t$ is an~asymptotically
finite-dimensional bounded semigroup then every $n$-dimensional
subspace $Y\subset X$ which complements~$X_0$ in~$X$
is {\it almost stabilizable}; i.e., the position of the space~$Y$
in the space $X$ varies under the action of the semigroup
successively slower:
$$
\forall t\quad
\sup\limits_{q<t}\angle(Y_T,Y_{T+q})\to 0
\quad
\text{as }T\to\infty.
$$
Moreover, the space $Y$ is not necessarily {\it stabilizable};
i.e., there may be no limit position~$Y_\infty$.
If such $Y_{\infty}$ exists then it is obviously
$\varphi_t$-invariant.
Also, we show ({\it Remark~3}) that the motion of an~almost
stabilizable but not stabilizable space~$Y$
in the space $X$ under the action of a~semigroup
cannot slow down too fast; and
the estimate from below for the change rate of the angle is the following:
$$
\sum\limits_{k=1}^{\infty}\angle(Y_{k+1},Y_k)=\infty.
$$

In {\bf Theorem~3} we prove that in the case when the semigroup
is {\it weakly almost periodic},
i.e., the closures of the orbits of vectors are compact in
the weak topology of~$X$ (for example,
in the case of a~bounded semigroup on a~reflexive~$X$)
the space $Y$ is stabilizable.

Using the criterion for invariance of a~finite-dimensional space
as a~proper space of the generator of a~semigroup
({\bf Lemma~6}), we construct {\bf Example~4}
which demonstrates that\linebreak boundedness of a~semigroup in Theorem~2
is essential even in the case
$\codim X_0=2$.

In the last section we prove {\bf Theorem~2$'$} which generalizes
Theorem~2 to the case of an~asymptotically finite-dimensional semigroup
with $\|\varphi_t\|=o(t)|_{t\to\infty}$
(these are exactly those semigroups
for which the summand $Q_t$ in representation~(1) is bounded).
This theorem\linebreak immediately implies that
if $\codim X_0=1$ (this case often emerges in applications)
then boundedness of a~semigroup in Theorem~2 is unessential.

Observe that analogs of Theorems~2 and 3  for
$C_0$-semigroups of operators are given in article~[1],
where they were proven by methods of nonstandard analysis.

The author is grateful to \`E.~Yu. Emel$'$yanov for
useful discussions.\vskip3mm

\centerline{\bf 1. Asymptotically Finite-Dimensional Semigroups}\vskip2mm

{\bf Lemma 1. }
Let $\varphi_t$ be a~bounded semigroup, $\|\varphi_t\|\leq C$.
Suppose that $v\in X$ and\linebreak $m(v)=\inf\limits_{t<\infty}\rho\{v_t, X_0\}=0$
{\rm (}i.e., among the vectors
 $\varphi_t(v)$,
there exist vectors arbitrarily close to the space~$X_0${\rm )}.
Then $v\in X_0$.
In particular, $X_0$ is closed in~$X$.

{\it Proof}.
Assume that $\varepsilon>0$. If $m(v)=0$ then
$|v_q-x|<\varepsilon$
for some $q\geq 0$
and $x\in X_0$.
Then $|v_{q+t}-x_t|<C\varepsilon$ for all $t<\infty$.
At the same time,
$|v_{q+t}|-|v_{q+t}-x_t|\leq|x_t|\underset{t\to\infty}\longrightarrow0$.
Therefore,
$\limsup\limits_{t\to\infty}\{|v_t|\}\leq C\varepsilon$.
The number $\varepsilon$ is arbitrary; therefore,
$\limsup\limits_{t\to\infty}\{|v_t|\}=0$ and $v\in X_0$.
$\square$

If  $\codim X_0<\infty$ then the lemma is also valid an~unbounded
semigroup~$\varphi_t$.
In~[2] this was proven for the semigroup of powers of an~operator
in a~complex space.
In the same article,
there is a~counterexample to the conclusion of Lemma~1
in the case when the space $X_0$ has infinite codimension
and the semigroup is unbounded.

We give an~example of a~semigroup with a~nonclosed~$X_0$.

{\bf Example 1  (V. V. Ivanov). }
The space $X_0$ of the unbounded discrete semigroup\linebreak
$\{T^n:l_2\to l_2\mid n\in\Bbb N\}$,
where $T(x_1,x_2,x_3,\dots)=(2x_2,2x_3,2x_4\dots)$
contains all finite sequences and therefore is dense in~$l_2$.
However, the reader can easily note that $X_0\neq l_2$.

{\bf Theorem ~1. }
Suppose that $\varphi_t:X\to X$ is an~asymptotically
countable-dimension semigroup.
Then the subspace $X_0$ is closed and $\codim X_0<\infty$.
If $\varphi_t$ is a~$C_0$-semigroup
{\rm (}i.e., for each $v\in X$ the function $t\mapsto v_t$
is continuous at zero{\rm )}
then $\varphi_t|_{X_0}:X_0\to X_0$ is a~bounded semigroup.

{\it Proof}.
Suppose that $\varphi_t$ is a~$C_0$-semigroup.

It is well known [3] that the passage to subspaces
of countable codimension preserves the property of being barrelled.
Thus, the uniform boundedness principle holds for $X_0$.
For every point $v\in X_0$ the set $\{v_t\mid t\geq 0\}$
is bounded; therefore, there is a~number $C<\infty$ such that
$\|\varphi_t|_{X_0}\|\leq C$
for each $t>0$.
The operators $\varphi_t$  on  $\Cl(X_0)$ are bounded by the same constant.
Lemma~1 applied to the restriction of the semigroup to
the space~$\Cl(X_0)$
implies that $X_0=\Cl(X_0)$. According to Baire category theorem
the complete quotient-space~$X/X_0$
cannot be countable-dimensional; it is only finite-dimensional.

If $\{\varphi_t\}$ is not a~$C_0$-semigroup
then, instead of the set $\{v_t\mid t\geq 0\}$,
we can consider the set
$\{v_t\mid t\geq t_0\}$ for some $t_0>0$.
It follows from the uniform boundedness principle that all
operators~$\varphi_t$, $t\geq t_0$,
are uniformly bounded on~$X_0$.
From the arguments similar to those in the proof of Lemma~1
we conclude that the space~$X_0$ is closed in this case as well.
The theorem is proven.

{\it Remark 1. }
For $C_0$-semigroups the space $X_0$ is a~Banach range, because
it is complete with respect to the norm
$\|x\|:=\sup\{|x_t|\mid t\geq 0\}\geq |x|$.
Thereby the complementation principle enables us to prove
closure of~
$X_0$,
using the assumption
$\codim X_0<\infty $.
In general, the presence of a~closed (algebraic) complement to~$X_0$ in~$X$
implies closure of~$X_0$.

It follows from Theorem~1 that the norm of the space
in representations of the form\linebreak
$X=X_0\oplus Y$ is equivalent to the norm of the direct product. Therefore,
Lemma~1 and Remark~1 yield

{\bf Corollary. }
Suppose that an~asymptotically finite-dimensional semigroup
is represented  by expression~{\it (1)} and
$y\in Y$. If $\liminf\limits_{t\to\infty}Q_t(y)=0$ then $y=0$.
If the semigroup $\varphi_t$ is bounded
then the semigroup $Q_t:Y\to Y$ is bounded as well
{\rm  (}the converse is false, see the last section of the article{\rm )}.

{\bf Lemma 2. }
Let $\varphi_t$ be a~bounded semigroup, $\|\varphi_t\|\leq C$.
Let $v\in X$. The function\linebreak
$v\mapsto m(v):X\to\Bbb{R}$ defined in Lemma~{\rm 1}
is continuous.

{\it Proof}.
Suppose that $x$ and $y\in X$. For each  $\delta>0$,
there is $t$ such that $|y_t|\leq m(y)+\delta$. Then
$$
m(x)\leq|x_t|=|y_t+(x-y)_t|\leq |y_t|+C|x-y|\leq m(y)+\delta+C|x-y|.
$$
The number $\delta$ is chosen arbitrarily;
therefore, $m(x)\leq m(y)+C|x-y|$.
Interchanging $x$ and $y$ in the above arguments, we obtain
$|m(x)-m(y)|\leq C|x-y|$. $\square$

{\it Remark 2. }
The boundedness condition for a~semigroup in Lemma~2 is essential.
We give one example: $X=\Bbb R^2$, $\varphi_t(y,z)=(y+tz,z)$.
Then $X_0=0$. The function $m:\Bbb R^2\to \Bbb R^2$
is strictly positive outside zero
but discontinuous at~$(1,0)\in \Bbb R^2$, since $m(1,0)=1$ and
$m(1,-\varepsilon)=\varepsilon$. This observation
enables us to construct
a~``counterexample'' to Theorem~2 for an~unbounded semigroup (Example~4).

The following convention enables us to avoid writing down
redundant constants in\linebreak inequalities.
We say that a~value $F$ {\it has order of a~value\/}~$H$
if there is a~constant $k\in \Bbb R$
such that $F\leq k\cdot H$ under the described conditions.

{\bf Lemma 3. }
Suppose that $\varphi_t$ is an~asymptotically finite-dimensional
bounded semigroup,
$Y\subset X$
is an~$n$-dimensional subspace such that $X_0\oplus Y=X$,
and $e^1,\dots, e^n$ is a~basis for the space~$Y$.
There is $k>0$ such that $|y_t|\geq k(|\beta_1|+\cdots+|\beta_n|)$
for each $t\geq 0$ and every vector
$y_t=\beta_1e^1_t+\cdots+\beta_n e^n_t$.

{\it Proof}.
We have
$$
\frac{|y_t|}{|\beta_1|+\cdots+|\beta_n|}=
\frac{|y_t|}{|y|}\cdot\frac{|y|}{|\beta_1|+\cdots+|\beta_n|}.
$$
Boundedness from below of the first factor follows
from Lemma~2 and the fact that the space~$Y$ is finite-dimensional;
boundedness from below of the second factor is obvious. $\square$

According to this lemma, the coefficients~$\beta_i$
in the decomposition
$y_t=\beta_1e^1_t+\cdots+\beta_n e^n_t$
cannot be too large if $|y_t|\leq 1$.

{\bf Corollary. }
Let $Z$ be an~$n$-dimensional subspace in~$X$.
The angle $\angle(Y_t, Z)$ has order of the maximal
distance from the vector $e^i_t$ to~$Z$, $i=1,\dots, n$.

{\bf Theorem ~2. }
Suppose that $\varphi_t$ is an~asymptotically finite-dimensional
bounded semigroup on the space~$X$ and
$Y\subset X$
is an~$n$-dimensional subspace such that $X_0\oplus Y=X$.
Then $Y$ is almost stabilizable; i.e.,
$\sup\limits_{s<t}\angle(Y_T,Y_{T+s})\to 0$
for each $t<\infty$
as $T\to\infty$.

{\it Proof}.
Represent the action of the semigroup
$\varphi:X_0\times Y\to X_0\times Y$ by means of~(1).
Note that $Q_t:Y\to Y$ is a~semigroup.

Let $e^1,\dots, e^n$ be a~basis for the space~$Y\subset X$.
In this basis, the mappings $Q_s:Y\to Y$ are given
by the matrix $(q_{ij})_s$ whose columns are
constituted by the coordinates of the projections of the vectors~
$e_s^i$ to the space~$Y$ parallel to the space~$X_0$:
\begin{equation}
\left(\begin{array}{c}
  e_{s}^1  \\
  \vdots\\
  e_{s}^n
\end{array}\right)=
\left(\begin{array}{ccc}
  q_{11}&\cdots& q_{1n} \\
  \cdots& \cdots &  \cdots\\
  q_{n1}&\cdots& q_{nn}
\end{array}\right)_s
\left(\begin{array}{c}
  e^1  \\
  \vdots\\
  e^n
\end{array}\right)
+ \left(\begin{array}{c}
  x^1(s)  \\
  \vdots\\
  x^n(s)
\end{array}\right),
\quad   x^i(s)=b_s(e^i)\in X_0.
                                                                \end{equation}

Applying the operator $\varphi_T$ to expression~(2) we obtain

\begin{equation}
\left(\begin{array}{c}
  e_{T+s}^1  \\
  \vdots\\
  e_{T+s}^n
\end{array}\right)
=
\left(\begin{array}{ccc}
  q_{11}&\cdots& q_{1n} \\
  \cdots& \cdots &  \cdots\\
  q_{n1}&\cdots& q_{nn}
\end{array}\right)_s
\left(\begin{array}{c}
  e_T^1  \\
  \vdots\\
  e_T^n
\end{array}\right)
+
\left(\begin{array}{c}
  x^1(s)_T  \\
  \vdots\\
  x^n(s)_T
\end{array}\right)
                                                                \end{equation}

for each $s\in[0,t]$.

The vectors of the first summand
on the right-hand side of~(3) lie in the space~$Y_T$. By the corollary
to Lemma~3, the angle $\angle(Y_{T+s}, Y_T)$ has order
$f(s,T):=\max\{|x^1(s)_T|,\dots,|x^n(s)_T|\}$.
However, all $x^i(s)$ lie in~$X_0$; therefore,
$f(s,T)\underset{T\to\infty}\longrightarrow0$.

Suppose that $\varphi$ is a~$C_0$-semigroup; i.e.,
the functions of the form $t\mapsto v_t$ are continuous at zero as well.
Then the set $\{x^i(s)\mid s\in[0,t]\}\subset X_0$
is compact, being the continuous image of the interval
$[0,t]$. By the uniform boundedness principle,
$\sup\{f(s,T)\mid s\in[0,t]\}\underset{T\to\infty}\longrightarrow0$.
Theorem~2 is proven for $C_0$-semigroups.

If $\varphi_t$ is not a~$C_0$-semigroup
then the function $x^i(s)$ may have discontinuity at zero;
therefore, the set $\{x(s)\mid s\in[0,t]\}$ may fail to be compact.
In this case we use a~somewhat artificial method:
as the initial space~$Y$ we consider the space~$Y$ already
translated by the action of the semigroup~$\varphi_t$;
i.e., the space~$Y_p$, $p>0$. Then for each $e^i\in Y_p$
the function $t\mapsto e^i_t$ is continuous at zero as well,
since there exist vectors~$u^i\in Y$ such that
$e^i=u^i_p$ and consequently $e^i_t=u^i_{p+t}$. Therefore, the functions
$x^i(s)= b_{s}(e^i)=b_{p+s}(u^i)$
are continuous for $s\geq 0$
rather than only for $s>0$.
Theorem~2 is proven completely.

{\it Remark 3. }
If the stabilization rate of the space~$Y$ is sufficiently high
then the space is stabilizable, i.e., tends to some
limit stable position~$Y_\infty$.
Indeed, the space $G(X,n)$ of $n$-dimensional subspaces
of the Banach space $X$ with the angular metric is complete.
Therefore, if, for example,
$Y_k\in G(X,n)$
is a~Cauchy sequence then it has a~limit $Y_\infty\in G(X,n)$.
It follows from~Theorem~2 that the oscillation of the function
$Y_t:t\to G(X,n)$ on the interval $[k,k+1]$
is small at large~$k$.
Thereby
$Y_\infty=\lim\limits_{t\to\infty}Y_t$.

In particular, if $Y$ is not stabilizable then the series
$\sum\limits_{k=1}^{\infty}\angle(Y_{k+1},Y_k)$
diverges. At the same time, it may happen that $Y$ is stabilizable,
but the series diverges. This follows from the fact that
the Cauchy condition is weaker than the condition of absolute
convergence of a~series.
We illustrate the last two conclusions.

{\bf Example 2. }
Let
$X=C[0,1]$ and $(\varphi_t f)(x)=x^tf(x)$~[1]. Here
$X_0=\{f\in C[0,1]\mid f(1)=0\}$ and $\codim X_0=1$.
There are no invariant complementing spaces. Hence, for every
function $f\in X$, if $f(1)\neq 0$ then the series $\sum\|f_{k+1}-f_k\|$
diverges.
We present computations for the function $f(x)\equiv 1$:
$$
\|f_{k+1}-f_k\|_X=\sup\limits_{x\in[0,1]} |x^{k+1}-x^k|=
\frac{1}{k+1}\left(\frac{k}{k+1}\right)^k\underset{k\to\infty}\to{\sim}
\frac{1}{ek},\ \ \ \
\sum\|f_{k+1}-f_k\|=\infty.
$$

{\bf Example 3. }
Let $X$
be the subspace of~$C[0,\infty)$ constituted by the functions having
a~limit at infinity and let $(\varphi_tf)(x)=f(x+t)$.
Then $X_0$ is a~space of functions tending to zero.
The space of constant functions is invariant.
Let $f(x)=1+\frac{\sin\pi x}{x}$.
The sequence $f_k(x)$ converges to unity uniformly; therefore, the space~$Y$
spanned by the vector $f\in X$ is stabilizable.
However, $\|f_{k+1}-f_k\|\sim \frac2k$
and the corresponding series diverges.\vskip3mm

\centerline {\bf 2. Analysis of the Evolution of Vectors in the Weak Topology}\vskip2mm

We now discuss some facts concerning the behavior of vectors
under the action of a~semigroup in the weak topology of the space~$X$.

Denote by the symbol $\Cl_\sigma$
the operator of taking the weak closure. Given a~number\linebreak
$0\leq r<\infty$ and a~vector $e\in X$, we put
$E_r=\{e_t\mid t\geq r\}$.
In particular, $E_0$ is the orbit of the vector $e$
under the action of the semigroup~$\varphi_t$.

{\bf Lemma 4. }
Let $\varphi_t$ be an~asymptotically finite-dimensional
bounded semigroup and let $e\notin X_0$.
Then $\Cl_\sigma(E_0)\cap X_0=\emptyset$.

{\it Proof}.
Suppose that $X=X_0\oplus Y$, $e\in Y$, and $e\neq 0$.
The continuous projection operator~$P:X=X_0\times Y\to Y$
is also continuous in the weak topology. At the same time,
the weak and strong topologies coincide on a~finite-dimensional~$Y$.
Using Lemma~1, we easily see that the projection of the orbit~$E_0$
of the vector~$e$ to~$Y$ is separated from zero; therefore,
$E_0\subset P^{-1}P(E_0)$ is separated from~
$X_0$
even in the weak topology. $\square$

{\bf Lemma 5. }
Let $e\in X$. The set $E_\infty=\bigcap\limits_{r<\infty}\Cl_\sigma(E_r)$
{\rm(}perhaps, empty{\rm)}
is $\varphi_t$-invariant; i.e.,\vskip-2mm
$E_{\infty+t}=E_\infty$
for each
$t$.

{\it Proof}.
The condition $z\in E_\infty$ means that
\begin{equation}
\forall \varepsilon>0\, \forall f_1,\dots, f_k\in X'\,
\forall r<\infty\  \exists T>r: |f_i(z-e_T)|<\varepsilon.
                                                                \end{equation}

Given a~functional $f\in X'$, we define the functional $f^t\in X'$
by the condition\linebreak
$f^t(x):=f(x_t)$. Then, applying condition~(4) to the functionals
$f^t_1,\dots, f^t_k\in X'$,
we find that there is an~arbitrarily large number~$T$ such that
$\bigl|f^t_i(z-e_{T})\bigr|<\varepsilon$
for each $i=1,\dots,k$. However,
$f^t_i(z-e_{T})=f_i(z_t-e_{T+t})$.
Thus, the vector $z_t$ satisfies condition~(4) and
$z_t\in E_{\infty}$. $\square$

{\bf Theorem ~3. }
Suppose that $\varphi_t$ is an~asymptotically finite-dimensional bounded
semigroup such that the orbit~$E_0$ of every vector $e\in Y$
is weakly precompact.
Then the space $Y$ is {\it stabilizable}; i.e.,
there is a~$\varphi_t$-invariant subspace~$Y_\infty$
such that $X= X_0\oplus Y_\infty$ and\linebreak
$\angle(Y_T, Y_\infty)\underset{T\to\infty}\longrightarrow0$.
In particular, if $X$ is reflexive then $Y$ is stabilizable.

The idea of the proof (the one-dimensional case):
we can approximate a~vector $z$ in the nonempty set
$ E_\infty=\bigcap\limits_{r<\infty}\Cl_\sigma(E_r)$
by convex combinations of vectors of the form~$e_{t_j}$
with arbitrarily large~$t_j$.
It follows from Theorem~2 that such combinations vary
arbitrarily slow. Hence, the vector
$z$
does not vary at all.
The multidimensional case is more complicated
only in the technical respect: we have to take care of the motion
of the space $Y$ with respect to itself.

{\it Proof}.
Consider the maximal set of vectors
$e^1=e$, $e^2=\varphi_{p_2}e^1,\dots, e^{s}=\varphi_{p_s}e^1$,  $s\leq n$,
whose projections to the space $Y$ are linearly independent.
It is easy to see that these projections constitute a~basis for
some $s$-dimensional subspace $Y^s\subset Y$
invariant under the action of the semigroup
$Q_t:Y\to Y$. We assume that $s=n$.
(In the general case the space~$Y$ is representable
as the direct sum of subspaces of the form~$Y^s$
and we apply the procedure described below to each of them.)

The set $\Cl_\sigma(E_0)$ is compact in the weak topology.
Hence, the intersection $E_\infty$ of the family of nested sets
$\bigcap\limits_{r<\infty}\Cl_\sigma(E_r)$ is nonempty.
Let $z^1\in E_\infty$. By~Lemma~4, $z^1\notin X_0$.

By Mazur's theorem, the weak closure of the set
lies in the closure of its convex hull.
Thus, it follows from the condition
``$z^1\in \Cl_\sigma(E_{T})\ \forall T<\infty$'' that,
whatever $\varepsilon>0$ might be, for every $T<\infty$,
there exist numbers
$\alpha_1,\dots,\alpha_m\geq 0$,
$\sum \alpha_k=1$ and vectors $e^1_{t_1},\dots,e^1_{t_m}$,
$t_j>T$, such that
$$
 z^1-\sum\limits_{k=1}^m\alpha_k e_{t_k}^1=\widetilde{\varepsilon}^1,
\quad |\widetilde\varepsilon^1|<\varepsilon.
$$

Applying the operator $\varphi_{p_i}$ to the last expression and denoting
$z^i=z^1_{p_i}$,
$i=1\dots,n$ (recall that $e^i=e^1_{p_i}$, $i=1,\dots, n$),
we obtain
\begin{equation}
z^i-\sum\limits_{k=1}^m\alpha_k e_{t_k}^i
=\widetilde{\varepsilon}^i,
\quad  |\widetilde{\varepsilon}^i|\leq C\varepsilon. \end{equation}\vskip-5mm

By~(3),\vskip-5mm
$$
e^i_{t_k+t}=\sum\limits_{j=1}^n q_{ij}e^j_{t_k}+x^i_{t_k},
$$
where $x^1,\dots,x^n\in X_0$.
Therefore,
\begin{equation}
\gathered
\sum\limits_{k=1}^m\alpha_k e_{t_k+t}^i=
\sum\limits_{k=1}^m\alpha_k\Biggl(\sum\limits_{j=1}^n
q_{ij}e_{t_k}^j+x^i_{t_k}\Biggr)=
\sum\limits_{j=1}^n q_{ij}\sum\limits_{k=1}^m\alpha_ke_{t_k}^j
+\sum\limits_{k=1}^m\alpha_kx^i_{t_k}
\\
=\sum\limits_{j=1}^n q_{ij}(z^j-\widetilde\varepsilon^j)
+\sum\limits_{k=1}^m\alpha_kx^i_{t_k}=
\sum\limits_{j=1}^n q_{ij}z^j-\sum\limits_{j=1}^n
q_{ij}\widetilde\varepsilon^j+\sum\limits_{k=1}^m\alpha_kx^i_{t_k}.
\endgathered
                                                                \end{equation}

Note that
$$
\Biggl|\sum\limits_{k=1}^m\alpha_kx^i_{t_k}\Biggr|
\leq \max\limits_{k=1,\dots, m}\bigl|x^i_{t_k}\bigr|.
$$
It follows from~(5) and (6) that
\begin{equation}
\Biggl|z_t^i-\sum\limits_{j=1}^n q_{ij}z^j\Biggr|
<\Biggl|\sum\limits_{j=1}^n q_{ij}\widetilde\varepsilon^j\Biggr|+
|\widetilde\varepsilon^i|+ \max\limits_{k=1,\dots, m} \bigl|x^i_{t_k}\bigr|.
                                                                \end{equation}

Since
$t_k>T$,
choosing $T$ sufficiently large,
we can make the whole right-hand side of inequality~(7)
be of order~$\varepsilon$.
However, the left-hand side of~(7) is independent of~$\varepsilon$;
therefore,
$$
z_t^i=\sum\limits_{j=1}^n q_{ij}z^j.
$$
Consequently, for each $t>0$, the vector $z_t^i$ lies in the linear
span of the vectors $z^1,\dots, z^n$;
hence, the space $Y_\infty$ spanned by the vectors $z^1,\dots, z^n$
is invariant.
The rest is obvious. The theorem is proven.

It follows from Lemma~5 that  $Y_\infty$ is
the linear span of the set~$E_\infty$.\vskip3mm

\centerline{\bf 3.   An~Infinitesimal Criterion for Invariance
and~Nonstabilizability}\vskip2mm

In this section we present an~infinitesimal criterion
for invariance of finite-dimensional subspaces and, using its,
construct an~asymptotically two-dimensional semigroup (Example~4)
having stable and nonstabilizable subspaces which complement~$X_0$ in~$X$.

For each semigroup $\varphi_t:X\to X$ denote by
$\varphi:X\to X$
the infinitesimal generator of~$\varphi_t$;
i.e., $v\in \dom \varphi$ if the limit\vskip-3mm
$$
\varphi(v)=\lim\limits_{t\to 0}\frac{\varphi_t(v)-v}{t}
$$
exists.
All properties of~$\varphi$ used below in the proofs
can be found, for example, in~[4].

{\bf Lemma 6 (a~criterion for invariance). }
Let $\varphi_t$ be a~semigroup.
All finite-dimensional $\varphi_t$-invariant subspaces of~$X$
are proper finite-dimensional subspaces
of the generator~$\varphi$ lying in~$\dom \varphi$.

{\it Proof}.
Suppose that $Y\subset \dom \varphi$ and $\varphi(Y)\subset Y$. For each
$y\in Y$
we have
$y_t=\sum\limits_{n=0}^\infty\frac{t^n\varphi^n(y)}{n!}\in Y$.
Therefore, $Y_t\subset Y$.
Conversely, suppose that $\dim Y<\infty$ and
$Y_t\subset Y$
for each $t<\infty$.
The restrictions $\psi_t=\varphi_t|_Y:Y\to Y$
constitute a~semigroup acting on~$Y$. The space~$Y$
is finite-dimensional; therefore, $\psi_t=e^{t\psi}$
and the infinitesimal operator $\psi:Y\to Y$
is defined everywhere.
But $\psi$ is a~restriction of the operator~$\varphi$. $\square$

Let $\alpha_t$ be a~semigroup on the space~$X$.
If $Q_t:B\to B$ is some semigroup with a~generator~$Q$
and $P:B\to X$ is a~continuous operator then the operator
$\varphi=\left(\begin{array}{cc}
\alpha & P \\
0& Q\end{array}\right)$
generates the semigroup
$\varphi_t:X\times B\to X\times B$
given by the formulas\vskip-1mm
\begin{equation}
\varphi_t
\left(\begin{array}{c}
  x \\
  b
\end{array}\right)=
\left(\begin{array}{cc}
  \alpha_t & \int\limits_0^t \alpha_s PQ_{t-s}d\,s \\
  0& Q_t
\end{array}\right)
\left(\begin{array}{c}
  x \\
  b\end{array}\right).
 \end{equation}

Suppose that
$B=\Bbb R^2$ and $Q_t:\Bbb R^2\to \Bbb R^2$ is defined as in Remark~2:\vskip-1mm
\begin{equation}
Q_t
\left(\begin{array}{c}
y  \\
z
\end{array}\right)
  =\left(\begin{array}{cc}
  1 & t \\
  0 & 1
\end{array}\right)
\left(\begin{array}{c}
y  \\
z\end{array}\right)
=\left(\begin{array}{c}
y+tz  \\
z\end{array}\right).
\end{equation}

Take $g\in X$. Define the mapping $P:B\to X$ by the formula\vskip-2mm
\begin{equation}
P(y,z)=y\cdot g.
\end{equation}

Define the semigroup $\varphi_t:X\times B\to X\times B$
by~(8).  The corresponding generators are as follows:\vskip-5mm
\begin{equation}
Q= \left(\begin{array}{cc}
  0 & 1 \\
  0 & 0
\end{array}\right),
\quad
\varphi\left(\begin{array}{c}
  f  \\   y\\ z\end{array}\right)=
\left(\begin{array}{c}
\alpha(f)+yg\\ z  \\ 0
\end{array}\right).
\end{equation}

Every space~$Y$ which complements~$X$ in~$X\times \Bbb R^2$
is the linear span of the vectors $(k,1,0)$ and $(l,0,1)$
for some $k, l\in X$.
Find out what are conditions on~$k$ and~$l$,
for $Y$ to be $\varphi_t$-invariant.

{\bf Lemma 7. }
Suppose that $\alpha_t:X\to X$ is a~semigroup, $g\in X$, and
$Q_t$ and $P(y,z)$
are defined by~{\rm (9)} and {\rm (10)}.
Let
$\varphi_t:X\times \Bbb R^2\to X\times \Bbb R^2$
be the semigroup defined by~{(8)}.  The vectors $u=(k,1,0)$ and $v=(l,0,1)$
generate a~$\varphi_t$-invariant subspace $Y$ if and only if
$k,l\in\dom\alpha$, $g=-\alpha(k)$, and $k=\alpha(l)$.

{\it Proof}.
Suppose that $u$, $v$ is a~basis for a~$\varphi_t$-invariant
space~$Y$. It follows from~Lemma~6 and~(11)
that $k,l\in\dom\alpha$ and
$\varphi(u)=(\alpha(k)+g, 0,0)=0$; i.e., $\alpha(k)+g=0$.
Now,\linebreak $\varphi(v)=(\alpha(l),1,0)=v$;
i.e., $\alpha(l)=k$. $\square$

{\bf Example 4. }
The semigroup of translations on $X=C_0(\Bbb R_+)$
is asymptotically zero-dimensional:
$$
X=X_0=C_0(\Bbb R_+)=\{f\in C(\Bbb R_+)\mid f(x)\to 0\},\quad
  (\alpha_tf)(x)=f(x+t),\quad  \alpha(f)=f'.
$$
Suppose that $g(x)\in X$ and the operators $Q_t$
and $P$ are the same as in~(9) and~(10).
Formula~(8) defines an~asymptotically two-dimensional semigroup
$\varphi_t:X\times \Bbb R^2\to X\times \Bbb R^2$:
\begin{equation}
\varphi_t
\left(\begin{array}{c}
f(x)\\ y\\ z
\end{array}\right)
=\left(\begin{array}{c}
  f(x+t)+\int\limits_0^tg(x+s)(y+(t-s)z)\,ds \\
  y+tz\\ z
\end{array}\right).
                                                                \end{equation}

It follows from Lemma~7 and the equality
$\alpha(f)=f'$
that the semigroup
$\varphi_t$
possesses an~invariant space which complements~$X$
if and only if the function $g(x)$ has the first and second
antiderivatives with the zero limit.

As~$g(x)$ consider the function $\frac{\sin(x)}{x}$.
Then the functions $k(x)=-\Si(x)+\frac\pi 2$
and $l(x)=x(\frac\pi 2-\Si(x))-\cos x$ satisfy the conditions of Lemma~7.
Examining the asymptotic behavior
of the integral sine function, we see that $k(x)$ and $l(x)$ vanish,
i.e., lie in~$X$. Thus, the corresponding semigroup has
a~two-dimensional invariant subspace which complements~$X$
in $X\times \Bbb R^2$.
Show that the subspace $Y=0\times\Bbb R^2\subset X\times \Bbb R^2$
is nonstabilizable.

{\bf Assertion. }
Let $g(x)=\frac{\sin(x)}{x}$. The space
$Y=(0\times \Bbb R^2)\subset X\times \Bbb R^2$
is not almost stabilizable under the action of semigroup~{\rm(12)}.

{\it Proof}.
Given a~vector $u\in Y_t$, denote by
$R(u)$ the vector joining~$u$ with the projection of~$u$
to the space $Y_{t+1}$ parallel to~$X$.
The key argument of Theorem~2 uses the fact that
$R(u)\underset{t\to\infty}\longrightarrow0$
for every
$u\in Y_t$;
moreover, convergence the is uniform; i.e.,
$$
\max\biggl\{\frac{|R(u)|}{|u|} \mid 0\neq u\in Y_t
\biggr\}\underset{t\to\infty}\longrightarrow 0.
$$
It turns out that the convergence may fail to be uniform
if the semigroup is unbounded.
In our example this is the case.
Indeed, suppose that
$v(t)=\varphi_t(0,-t,1)$.
Inserting the vector~$v(t)$ in~(12), we see that
$$
v(t)=
\left(\begin{array}{c}
-\int\limits_0^t\frac{\sin (x+s)}{x+s}s\,ds \\
  0\\1
\end{array}\right),
\quad
R(v(t))=v(t)-v(t+1)
=\left(\begin{array}{c}
\int\limits_t^{t+1}\frac{\sin (x+s)}{x+s}s\,ds \\
  0\\0
\end{array}\right).
$$

The value
$$
\int\limits_a^b\frac{\sin (x+s)}{x+s}s\,ds
=\cos|^{a+x}_{b+x} -x \Si|^{a+x}_{b+x}
$$
is bounded for all $x, a,b\geq 0$,
since
$\Si(p)\sim \frac{\pi}{2}-\frac{\cos p}{p}$ as $p\to\infty$.
Therefore,
$\frac{|R(v(t))|}{|v(t)|}\not\to0$
as $t\to\infty$, since
$$
|R(v(t))|=\left\|\int\limits_t^{t+1}\frac{\sin (x+s)}{x+s}s\,ds\right\|_X=
\max_{x\geq 0}\left|\int\limits_t^{t+1}\frac{\sin (x+s)}{x+s}s\,ds\right|
\overset{x:=0}\to{\geq}
|\cos|^t_{t+1}|\not\to 0.
$$

Thus, the distance from the vector
$v(t)=\varphi_t(0,-t,1)$ to its projection $v(t+1)$ to the
space $Y_{t+1}$ is sufficiently large.

It follows from~(12) that all elements of the space~$Y_t$ have the form
\begin{equation}
\left(\begin{array}{c}
h^t_{a,b}(x)\\a\\b
\end{array}\right)=
\left(\begin{array}{c}
\int\limits_0^t\frac{\sin (x+s)}{x+s}\cdot (a-sb)\,ds \\
  a\\b
\end{array}\right).
                                                                \end{equation}

We have shown that there is a~positive constant $K$ such that
there exists an~arbitrarily large number $t<\infty$ such that
the distance from the vector $v(t)=\bigl(h^t_{0,1},0,1\bigr)\in Y_t$
to its projection~$\bigl(h^{t+1}_{0,1},0,1\bigr)$
to $Y_{t+1}$ is greater than~$K$.
It remains to observe that for such $t$ the distance
from $v(t)$ to all other vectors
$\bigl(h^{t+1}_{y,z},y,z\bigr)\in Y_{t+1}$:
\begin{equation}
\rho_X\bigl\{v(t)-\bigl(h^{t+1}_{y,z},y,z\bigr)\bigr\}=
\max_{0<x<\infty}\bigl\{\bigl|h^t_{0,1}-h^{t+1}_{y,z}\bigr|\bigr\}
+\sqrt{y^2+(z-1)^2},
                                                                \end{equation}

is also bounded from below for all $(y,z)\in\Bbb R^2$.
This follows from the fact that, for $(y,z)$ close
to~$(0,1)$, the first summand in~(14)
cannot become small immediately
(the necessary estimate is simple and left to the reader),
while the second summand becomes more essential
as $(y,z)$ recedes further from~$(0,1)$.

Symmetric arguments show that the distance from the vector
$\bigl(h^{t+1}_{0,1},0,1\bigr)$ to the space~$Y_t$
is bounded from below as well.
Thus, $\angle(Y_t,Y_{t+1})\not\to 0$ as $t\to\infty$.\vskip3mm

\centerline{\bf 4.   Stabilizability in Slowly Increasing Semigroups}\vskip2mm

In this section all semigroup are supposed to
be asymptotically finite-dimensional.

The boundedness condition for the semigroup $Q_t:Y\to Y$
in representation~(1) is inde\-pen\-dent of the choice of the subspace~$Y$
which complements~$X_0$ in~$X$.
We call such semigroups {\it slowly increasing}.
Analysing the upper right entry of the matrix in representation~(8),
we can easily derive estimates which show that
the slow growth condition for the semigroup~$\varphi_t$
is equivalent to the condition
$\|\varphi_t\|=o(t)|_{t\to\infty}$.

{\bf Example 5. }
Consider the semigroup defined by formula~(8),
where $X_0=C_0(\Bbb R_+)$
is the same as in Example~4, $Y=\Bbb R$, and $P(y)=g\cdot y$, where
$g\in X_0$:
\begin{equation}
\varphi_t
\left(\begin{array}{c}
f(x) \\
y
\end{array}\right)=
\left(\begin{array}{c}
f(x+t)+y\cdot\int\limits_0^tg(x+s)\,ds\\
y
\end{array}\right).
                                                                \end{equation}

The function $g(x)$ can be such that its integral
may attain arbitrarily large values.
At the same time, the growth rate of~
$\|\varphi_t\|$ is determined by the growth rate of the function
$\int\limits_0^t g(t)\,dt$.
Therefore, semigroup~(15) may fail to be bounded.
At the same time, $\|\varphi_t\|=o(t)$ as
$t\to\infty$, since $g(x)\to 0$.
For example, if  $g(x)=\frac{1}{x+1}$ then, putting
$f\equiv 0$ and $y=1$ in~(15), we obtain
$$
\varphi_t(0,1)=\biggl(\ln\frac{x+1+t}{x+1},1\biggr),
\quad
\|\varphi_t\|\geq\sim \sup\biggl\{\ln\frac{x+1+t}{x+1}
\mid x\geq 0\biggr\}\sim\ln t.
$$

Thus, the class of slowly increasing semigroups is broader
than the class of bounded\linebreak semigroups.
Show that the conclusion of Theorem~2 holds
in this class as well.

{\bf Theorem ~2$'$. }
The conclusion of Theorem~{\rm 2} is valid
for all slowly increasing semigroups.

{\it Proof}.
Suppose that $X=X_0\oplus Y$, $\dim Y<\infty$,
and
$\varphi_t=
\left(\begin{array}{cc}
\alpha_t &  b_t \\
  0 & Q_t
\end{array}\right) $.

First consider the case $\dim Y=1$. On the one hand,
this case presents the main idea of the proof in the general case
and, one the other hand, is rather simple;
in this case the semigroup $Q_t:Y\to Y$
is $\exp(ct)$ for some $c\in\Bbb R$.
By the conditions of the theorem, $c\leq 0$.
According to the corollary to Theorem~1, $c=0$.
Assume that $y\in Y$ and $t>0$. Then $y_t=x(t)+y$,
where $x(t)\in X_0$. Thus,
$|y_t|\geq |y|$.
At the same time,\vskip-3mm
$$
|\varphi_T(y_t)-\varphi_T(y)|=|\varphi_T(y_t-y)|=
|\varphi_T(x(t))|\underset{T\to\infty}\longrightarrow  0.
$$\vskip-1mm
This shows that the angle between the straight line $Y_t$ and
$Y_{T+t}$ vanishes at large~$T$.

Turning to the general case, observe that
if the semigroup $Q_t:Y\to Y$ is bounded then
it is also bounded from below, moreover, uniformly.
We can derive this from the corollary to Theorem~1 and
Lemma~2 applied to the {\it finite-dimensional\/} semigroup $Q_t:Y\to Y$
itself.
Thus, there is $k<\infty$ such that
$|y|\leq k|Q_t(y)|$
for all $y\in Y$ and $t>0$. Hence,
 $|y|\leq k|\varphi_t(y)|$.

Assume that $t>0$. The ball $B\subset Y$ of radius $k$ is compact,
while its image~$B_t$ is also compact. Then the set
$A:=X_0\cap (B_t-B)=\{u-v\in X_0 \mid u\in B_t,\ v\in B\}$
is compact as well.
In the rest of the proof, the set~$A$ plays
the same role as the point~$x(t)$ in the proof of the one-dimensional case.

Suppose that $z\!\in\! Y_T$ and~$|z|\!=\!1$.
Consider a~vector
$y\in Y$
such that  {$z=y_T$.
Then $|y|\leq k$;} i.e., $y\!\in\! B$.
There is $x\!\in\! X_0$ such that $y+x\!\in\! Y_t$. Then $x\!\in\! A$.
At the same time,
$z+x_T\!=\!(y+x)_T\in Y_{T+t}$.
Therefore,\vskip-8mm
$$
\rho(z, Y_{T+t})\leq x_T\leq|A_T|:=\sup\{|x|\mid x\in A_T\}.
$$\vskip-1mm
The number $|A_T|$ is independent of the choice of~$y_T$; therefore,
$\angle(Y_T,Y_{T+t})\leq |A_T|$ by the definition of the angle.
But the set $A$ is compact and lies in~$X_0$;
therefore, $|A_T|\underset{T\to\infty}\longrightarrow 0$.
Consequently, $\angle(Y_T,Y_{T+t})\underset{T\to\infty}\longrightarrow 0$.
We are left with applying once again the uniform boundedness principle,
arguing as at the end of the proof of Theorem~2.  $\square$

Example~4 shows that the slow growth condition
for the semigroup $Q_t$ in Theorem~2$'$
is essential even in the case $\codim X_0=2$.
However, we can require nothing for asymptotically {\it one-dimensional\/}
semigroups~$\varphi_t$.
Indeed, as was observed, the mapping $Q_t:Y\to Y$
in repre\-sen\-ta\-tions~(1) of such semigroup
is the multiplication by the number~$e^{ct}$.
The operators $\psi_t:=e^{-ct}\varphi_t$ as well constitute a~semigroup
which is ``homothetic'' to the original one and is slowly increasing.
It follows from Theorem~2$'$ that the subspace~$Y$
is almost stabilizable in the semigroup~$\psi_t$
and hence in the semigroup~$\varphi_t$.
Moreover, if the semigroup $\psi_t$ is unbounded
(for example, as in Example~5)
then the angle between the straight line~$Y_t$ and the space $X_0$ vanishes.

In conclusion, note that if an~asymptotically
finite-dimensional semigroup is slowly increasing
but unbounded then there are no stable subspaces which
complements~$X_0$.
Indeed, the {\it direct} product of bounded semigroups is bounded.
Therefore the representation~(1) for such semigroup cannot be diagonal.\vskip-2mm

\centerline{\bf References}

1. Emel'yanov \`{E}.Yu., Some conditions for a $C_0$-semigroup to be asymptotically
finite-dimensional. Siberian Mathematical Journal, 44, \No 5, p.793--796 (2003).

2. Emel'yanov \`{E}.Yu. and Wollf M.P.H., Quasi constricted linear operators on Banach spaces.
Studia Math., 144, \No 2, p. 169--179 (2001).

3. Levin M and Saxon S., Every countable-codimensional subspace of a barreled space is barreled.
PAMS, 29, \No 1, p.91--96 (1971).

4. Hille E. and Phillips R.S., Functional Analysis and Semi-Groups. AMS, colloquium publ., v.XXXI.
Providence, 1957.

\enddocument